\documentclass[a4paper,10pt]{article}

\usepackage{setspace}
\usepackage[bottom]{footmisc}
\usepackage{color}
\usepackage{subfig}
\usepackage{hyperref}
\usepackage{multirow}
\usepackage{cite}
\usepackage{amsmath} 
\usepackage{mathtools}
\usepackage{cleveref}
\usepackage{graphicx}
\usepackage{epsfig}
\usepackage{epstopdf}
\usepackage{float}
\usepackage{titling}
\usepackage[margin=1.0in]{geometry}
\usepackage{etoolbox}
 
\title{Shock capturing with discontinuous Galerkin Method using Overset grids for two-dimensional Euler equations}
\makeatletter
\renewcommand\@date{{%
  \vspace{-\baselineskip}%
  \large\centering
  \begin{tabular}{@{}c@{}}
    S R Siva Prasad Kochi\textsuperscript{1} \\
    \normalsize siva.ksr@gmail.com
  \end{tabular}%
  \quad and\quad
  \begin{tabular}{@{}c@{}}
    M Ramakrishna\textsuperscript{2} \\
    \normalsize krishna@ae.iitm.ac.in
  \end{tabular}

  \bigskip

  \textsuperscript{1}Doctoral Candidate, Dept. of Aerospace Engg., IIT Madras.\par
  \textsuperscript{2}Professor, Dept. of Aerospace Engg., IIT Madras.

  \bigskip

  \today
}}
\makeatother
\begin{document}

\maketitle

\begin{abstract}
 A new procedure to capture the shocks has been proposed and is demonstrated for the solutions of two-dimensional Euler equations using discontinuous Galerkin method and overset grids. A discontinuous Galerkin solver using a coarse grid provides the troubled cell data that is used to determine the location of the shock. An overset grid aligned to the shock is constructed based on this information. The solver is run again with the coarse grid solution as the initial condition with limiting occurring only in the overset grid to give a solution with the shock aligned to the grid line. Results for supersonic flow over a ramp, shock reflecting off a flat plate and the supersonic flow over a circular cylinder are presented. For the flow over a circular cylinder, the results obtained are validated using an existing analytical method for calculating the shock offset distance.
 
 {{\bf Keywords:} discontinuous Galerkin method, shocks, overset grids}
\end{abstract}

\section{Introduction}

In this paper, we propose a procedure for capturing shocks and demonstrate it for the solutions of two-dimensional Euler equations using discontinuous Galerkin Method (DGM) and overset grids. Flows featuring shocks can be modeled by either shock capturing or shock fitting. Shock fitting \cite{salas1} is done by locating and tracking the motion of the shocks which are treated as boundaries between regions where smooth solutions exist. Shock fitting algorithms are in general complex where as shock capturing algorithms are quite simple and hence they have enjoyed a lot of popularity in recent years \cite{bpm}.
\\
\\
\noindent For higher order methods like DGM, a lot of shock capturing methods are proposed in literature. The main ingredients of such methods are shock detection and shock stabilization which have been studied extensively with regard to Finite Difference (FD) and Finite Volume (FV) methods. For DG methods, a number of shock detectors have been developed whose performance is compared by Qiu and Shu \cite{qs2}. The main approach to stabilize shocks is to increase mesh resolution in the detected regions. This is not always feasible. Mesh refinement is usually supplemented or completely replaced by limiting \cite{qs1} or artificial viscosity techniques \cite{pp}. A comparative study of such techniques has been done by Jian Yu and Hesthaven \cite{yh}. Each of these techniques have their pros and cons and selection of a particular technique is a matter of choice. We have used limiters in our shock capturing procedure.
\\
\\
\noindent Weighted Essentially Non-Oscillatory (WENO) limiters \cite{qs1} are the preferred limiters for DGM because they maintain the order of accuracy in cells where the solution is smooth. Accuracy of limiting is increased if they are used along with mesh refinement in the vicinity of the shock \cite{gq}. All these shock capturing procedures agree on the fact that if the grid lines are aligned to the shock, we can capture it more accurately. We propose the following procedure for shock capturing. We locate the shock using a good shock detector and use an overset grid so that the grid lines are aligned to the shock to capture it. We use the limiter only on the overset grid. This procedure can be called as shock capturing using shock fitting as we locate the shock (as in shock fitting) and change the grid so that the grid lines are aligned to the shock and capture it accurately.
\\
\\
\noindent Overset grids were first applied to solving the Euler equations by Benek et al.\cite{bsd}. They have used a body fitted grid to impose boundary conditions for a curved geometry. They produced the solution on a system of grids that communicate through exchange of boundary data (called artificial boundaries) in overlapping regions. The arbitrary overlapping of grids allows the mesh generator to focus on resolving individual components of the geometry independently. This provides great flexibility for the mesh generation process to focus resources on regions of interest. We use overset grids for shock capturing.
\\
\\
\noindent Our proposed method to use overset grids is to run the solver on a coarse grid first and use a good troubled cell indicator to find the location of the shock. Then, we construct an overset grid so that the grid lines are parallel to the shock. We run the solver again and the limiting process is used only in the overset grid. This will enable us to capture the shock while applying the limiter only where it is required.
\\
\\

\noindent The paper is organized as follows. We describe the formulation of the discontinuous Galerkin method used for all our results in \cref{sec:formulation}, the proposed procedure for shock capturing using overset grids is described in \cref{sec:oversetExpl}, and the validation of the overset grid solver is described in \cref{sec:validationOverset} and the results are described in \cref{sec:results} and we conclude the paper in \cref{sec:conc}.

\section{Formulation of discontinuous Galerkin Method}\label{sec:formulation}

\noindent Consider the Euler equations in conservative form as given by

\begin{equation}\label{2dEulerEquations}
\frac{\partial \textbf{Q}}{\partial t} + \frac{\partial \textbf{F(Q)}}{\partial x} + \frac{\partial \textbf{G(Q)}}{\partial y} = 0 \quad \text{in the domain} \quad \Omega
\end{equation}
\noindent where $\textbf{Q} = (\rho, \rho u, \rho v, E)^{T}$, $\textbf{F(Q)}=u\textbf{Q} + (0, p, 0, pu)^{T}$ and $\textbf{G(Q)}=v\textbf{Q} + (0, 0, p, pv)^{T}$ with $p = (\gamma -1)(E-\frac{1}{2}\rho (u^{2}+v^{2}))$ and $\gamma = 1.4$. Here, $\rho$ is the density, $(u,v)$ is the velocity, $E$ is the total energy and $p$ is the pressure. We approximate the domain $\Omega$ by $K$ non overlapping elements given by $\Omega_{k}$. 
\\
\\
We look at solving \eqref{2dEulerEquations} using the discontinuous Galerkin method. We approximate the local solution as a polynomial of order $N$ which is given by:

\begin{equation}\label{nodalApprox}
 Q_{h}(r,s) = \sum_{i=0}^{N} \sum_{j=0}^{N} Q^{ij} \Phi_{ij}(r,s) = \sum_{i=0}^{N} \sum_{j=0}^{N} Q^{ij} P_{i}(r)P_{j}(s)
\end{equation}

\noindent where $r$ and $s$ are the local coordinates. The polynomial basis is given by $\Phi_{ij}(r,s) = P_{i}(r)P_{j}(s)$ where $P_{i}(r)$ and $P_{j}(s)$ are the one-dimensional Lagrange interpolation polynomials of order $N$ at the appropriate one-dimensional Gauss Legendre quadrature points. Now, using $\Phi_{ij}(r,s)$ as the test function, the weak form of the equation \eqref{2dEulerEquations} is obtained as

\begin{equation}\label{weakFormScheme}
 \sum_{i=0}^{N} \sum_{j=0}^{N} \frac{\partial Q^{ij}}{\partial t} \int_{\Omega_{k}} \Phi_{pq} \Phi_{ij} d\Omega + \int_{\partial \Omega_{k}} \hat{F} \Phi_{pq} ds - \int_{\Omega_{k}} \vec{F}.\nabla \Phi_{pq} d\Omega = 0
\end{equation}

\noindent where $\partial \Omega_{k}$ is the boundary of $\Omega_{k}$, $\vec{F} = (\textbf{F(Q)},\textbf{G(Q)})$ and $\hat{F} = \bar{F^{*}}.\hat{n}$ where $\bar{F^{*}}$ is the monotone numerical flux at the interface which is calculated using an exact or approximate Riemann solver and $\hat{n}$ is the unit outward normal. This is termed to be $\mathbf{P}^{N}$ based discontinuous Galerkin method.
\\
\\
\noindent Equation \eqref{weakFormScheme} is integrated using an appropriate Gauss Legendre quadrature and is discretized in time by using the fifth order Runge-Kutta time discretization given in \cite{butcher} unless otherwise specified. To control spurious oscillations which occur near discontinuities, a limiter is used with a troubled cell indicator. We have used the modified KXRCF troubled cell indicator and the WENO limiter proposed in \cite{lpr} for all our calculations.

\section{Using overset grids for shock capturing}\label{sec:oversetExpl}
\begin{figure}[htbp]
\begin{center}
\includegraphics[scale=0.8]{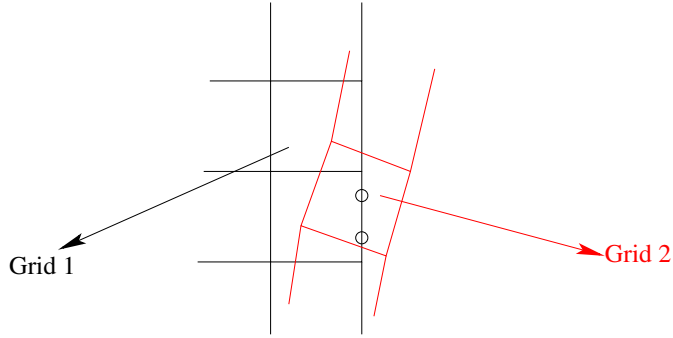}
\caption{Two overlapping grids (Grid 1 in black and Grid 2 in Red) with artificial boundary of Grid 1 containing Gauss Quadrature points}
\label{fig:OversetGridExpl}
\end{center}
\end{figure}

\noindent For a solver using overset grid, we have multiple grids which are overlapped with each other. Boundaries of two overlapping grids named Grid 1 (black) and Grid 2 (red) are shown in Figure \ref{fig:OversetGridExpl}. These boundaries are called artificial boundaries. Inter-grid communication happens through the artificial boundaries.
\\
\\

\begin{figure}[htbp]
  \centering
  \subfloat[Troubled cells (black) for Mach 3.0 flow over a $10^{\circ}$ ramp]{\label{fig1:rampM3p0troubledCells}\includegraphics[width=0.9\textwidth]{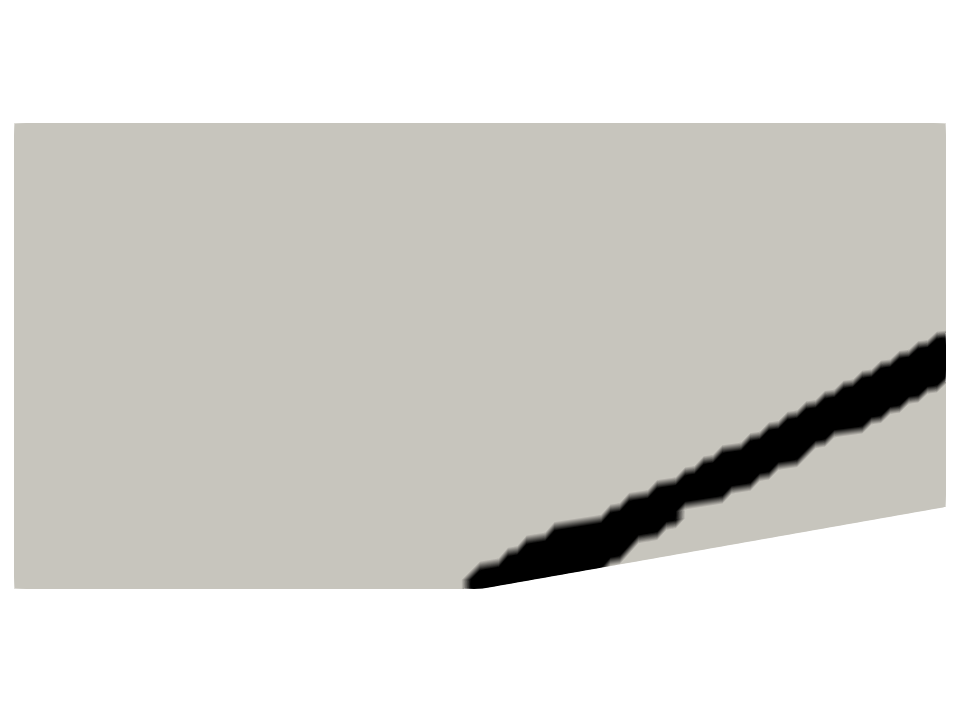}}\hfill
  \subfloat[Troubled cells (black) for Mach 3.0 flow over a circular cylinder]{\label{fig2:CircCylM3p0troubledCells}\includegraphics[width=0.9\textwidth]{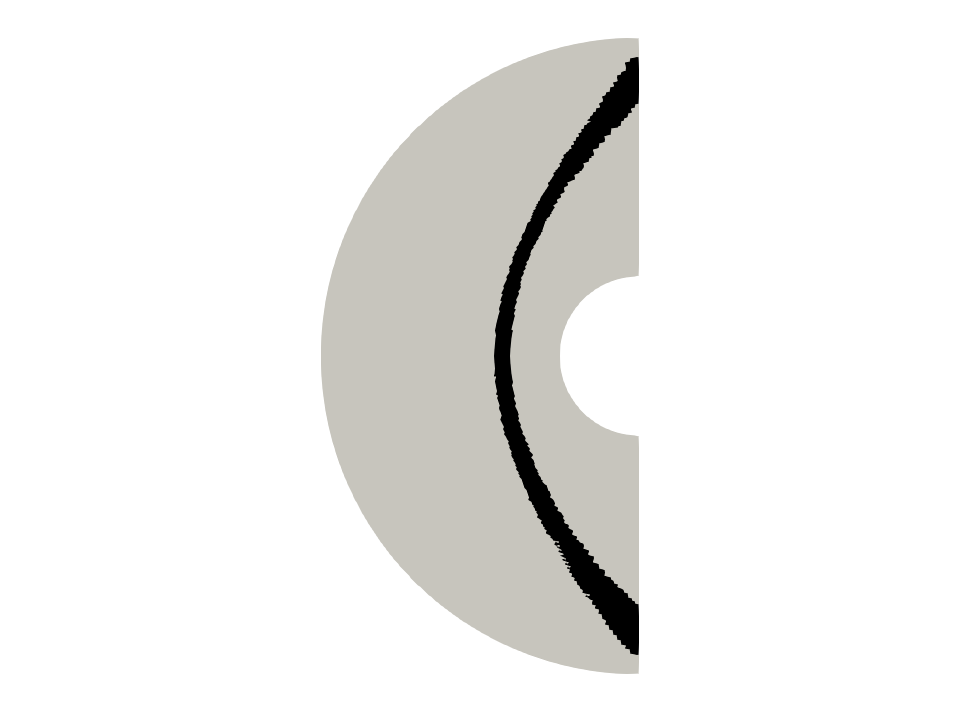}}\hfill
  \caption{Troubled Cells (black) using the modified KXRCF troubled cell indicator proposed in \cite{lpr}}
  \label{fig:troubledCells}
\end{figure}

\noindent For an artificial boundary, the exterior conservative variables $Q^{+}$ (where $Q$ is the solution variable)  must be provided by one, or multiple, elements from overlapping meshes for the calculation of numerical flux. We will explain how to obtain them now. As seen in Figure \ref{fig:OversetGridExpl}, one of the artificial boundary faces of Grid 1 is seeded with Gauss Quadrature (GQ) points as is required by the order of accuracy of the solver. We first find the Cartesian coordinates of these GQ points. These Cartesian coordinates are used to obtain the cell local coordinates in Grid 2. A search algorithm (K-d tree) is used to determine which GQ points are located in each cell of Grid 2. Then, we can find the cell local coordinates using Newton's method as given by \cite{gbot}. These coordinates are stored before hand so as not to calculate them at every time step. Using these cell local coordinates, we can find the required $Q^{+}$ values at each of the GQ points. This will be used as the boundary condition for Grid 1. We repeat the same procedure for Grid 2. This procedure is repeated for all overlapping grids at each time step.
\\
\\
\noindent We use overset grids to capture the stationary shocks that occur in the solution of Euler equations. The proposed step by step procedure is outlined below:
\\
\\
\noindent \textbf{Step 1:} Run the solver on a coarse grid with a given troubled cell indicator and limiter and obtain the solution.
\\
\\
\noindent \textbf{Step 2:} Look at the troubled cells to locate the discontinuities (shocks) that occur in the solution.
\\
\\
An example of troubled cell distribution for Mach 3.0 flow over a $10^{\circ}$ ramp and Mach 3.0 flow over a circular cylinder are shown in Figures \ref{fig1:rampM3p0troubledCells} and \ref{fig2:CircCylM3p0troubledCells} respectively. From these troubled cells, we get an idea of the location of the shock.
\\
\\

\begin{figure}[htbp]
\begin{center}
\includegraphics[scale=0.4]{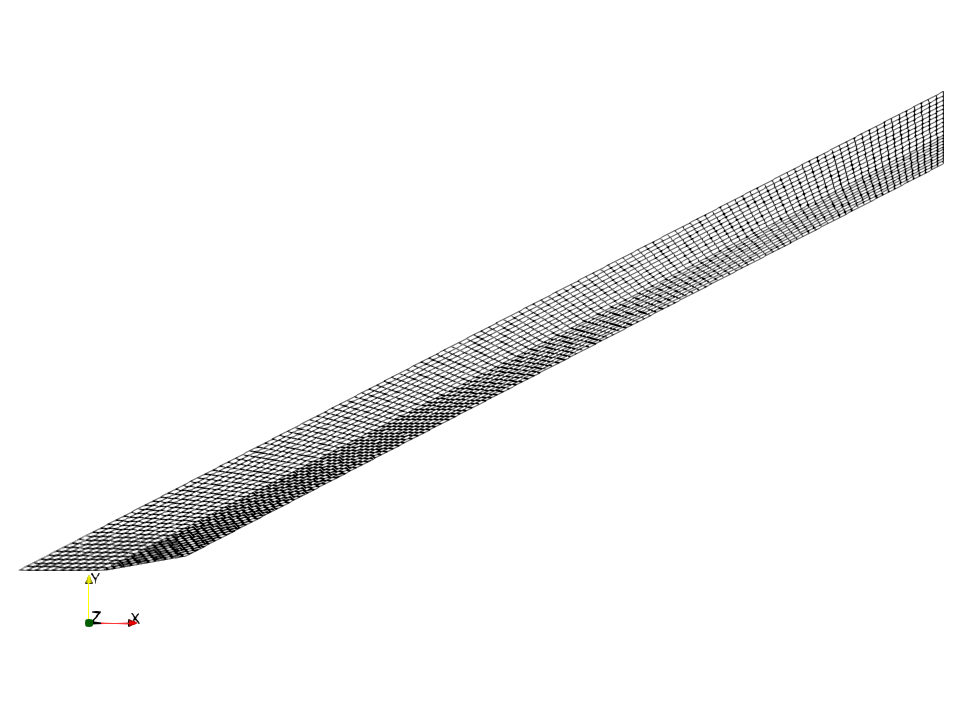}
\caption{Overset grid for Mach 3.0 flow over a $10^{\circ}$ ramp.}
\label{fig:rampOversetGrid}
\end{center}
\end{figure}

\begin{figure}[htbp]
  \centering
  \subfloat[Overset grid for Mach 3.0 flow over a circular cylinder]{\label{fig2:CircCylM3p0OversetGrid1}\includegraphics[width=0.45\textwidth]{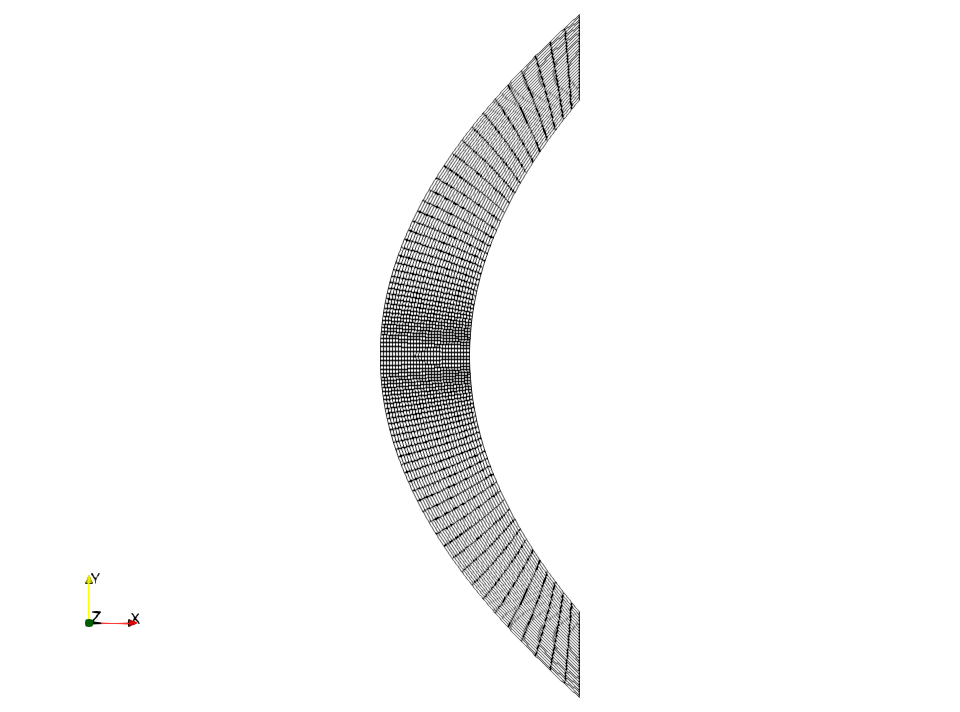}}\hfill
  \subfloat[Overset grid for Mach 3.0 flow over a circular cylinder with reduced width using procedure outlined in Step 4]{\label{fig3:CircCylM3p0OversetGrid2}\includegraphics[width=0.45\textwidth]{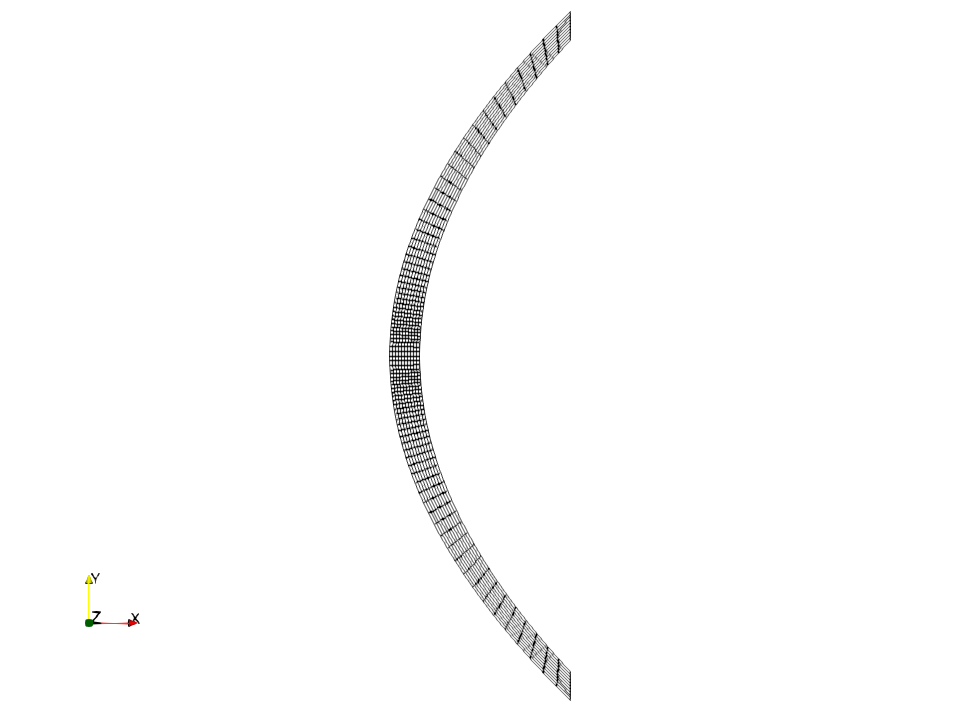}}\hfill
  \caption{Overset grid used for flow over circular cylinder}
  \label{fig:OversetGrid}
\end{figure}

\noindent \textbf{Step 3:} Construct an overset grid which encompasses all the troubled cells which is also confirming to the body. 
\\
\\
For the flow over a ramp, we know that the shock is straight and the overset grid is constructed to capture the shock. An example overset grid is shown in Figure \ref{fig:rampOversetGrid}. For the flow over a circular cylinder, the shock is curved. In this case, we use the cell centers of the first troubled cell that the flow encounters in the solution and interpolate those points using an appropriate b-spline approximation. This gives an approximate shape of the shock. We use this polynomial and construct a grid by using parallel offset polynomials so that the grid lines are parallel to each other. An example grid constructed in this way is shown in Figure \ref{fig2:CircCylM3p0OversetGrid1}.
\\
\\
%

\noindent \textbf{Step 4:} Using this overset grid, we run the solver again with the coarse grid solution as the initial condition. While running the solver, we need to use the troubled cell indicator and the limiter only inside the overset grid. We also use a high resolution numerical flux in the overset grid to capture the shock accurately. We have used the SLAU2 \cite{ks3} numerical flux in the overset grid and the Lax-Friedrichs flux elsewhere. After a few timesteps, we can look at the troubled cell profile again and reduce the width of the overset grid even further.
\\
\\
An example of such overset grid with reduced width is shown in Figure \ref{fig3:CircCylM3p0OversetGrid2}. While reducing the width, we have to take care that the coarse grid does not overlap the fine grid at the location of the shock. In this way, we can capture the shock accurately. If we want, we can use another overset grid over the shock to resolve solution near the shock even further. This completes the procedure.

\section{Validation of overset grid solver using Isentropic Vortex Problem}\label{sec:validationOverset}

\begin{figure}[htbp]
\begin{center}
\includegraphics[scale=0.4]{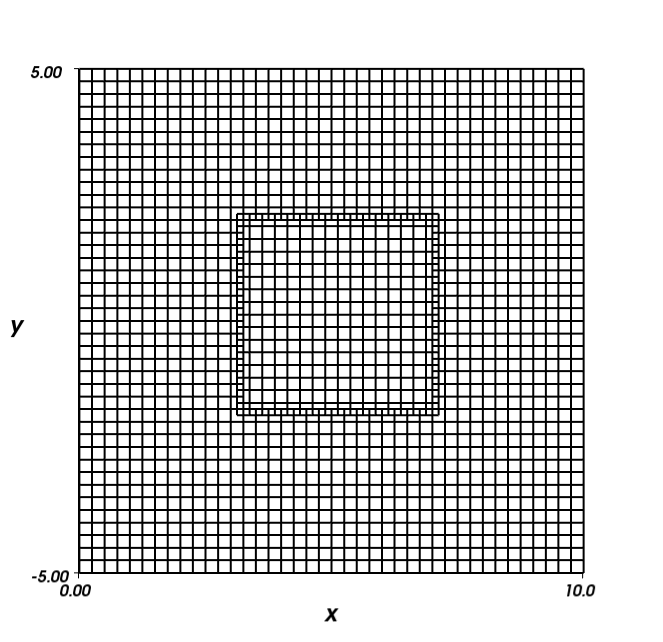}
\caption{Overset Grid for Isentropic Vortex Problem in the domain $[0,10]\times[-5,5]$ with $40$ by $40$ elements used for the validation of the overset grid solver.}
\label{fig:IsenVortexOversetGrid}
\end{center}
\end{figure}

\noindent In this section, we look at validating the overset grid solver to ensure that the data transfer that occurs in the overset grid does not decrease the order of accuracy of the scheme. We do this by using the Isentropic Euler Vortex problem \cite{shu1} which has a smooth solution. We run the solver on a given domain with uniform grid and also on a domain which contains an overset grid of the same size as the original grid as shown in Figure \ref{fig:IsenVortexOversetGrid}. We calculate the order of accuracy of the scheme for both grids to show that the effect of data transfer that occurs in the overset grid is negligible. We solve the two dimensional Euler equations given by \eqref{2dEulerEquations} in the domain $[0,10]\times[-5,5]$ for the Isentropic Euler Vortex problem. The exact solution is given by: \\ $\rho = \left(1 -  \left(\frac{\gamma - 1}{16\gamma \pi^{2}}\right)\beta^{2} e^{2(1-r^{2})}\right)^{\frac{1}{\gamma-1}}$, $u = 1 - \beta e^{(1-r^{2})} \frac{y-y_{0}}{2\pi}$, $v = \beta e^{(1-r^{2})} \frac{x-x_{0}-t}{2\pi}$, and $p = \rho^{\gamma}$, where $r=\sqrt{(x-x_{0}-t)^{2}+(y-y_{0})^{2}}$, $x_{0}=5$, $y_{0}=0$, $\beta=5$ and $\gamma = 1.4$. We initialize with the exact solution at $t=0$ and use periodic boundary conditions at the edges of the domain. We run the solver with the same grid size for a normal grid and an overset grid for grid sizes of $1/20$, $1/40$, $1/80$, and $1/160$ for various orders. An example overset grid used is shown in Figure \ref{fig:IsenVortexOversetGrid} for a grid of $40$ by $40$ elements. The errors in density and numerical orders of accuracy are calculated at $t=10.0$ (one period) for the original grid as well as the overset grid and are presented in Table \ref{table:1}. While calculating the solution, we have made sure that the temporal and spatial orders of accuracy are the same by using a corresponding Runge-Kutta time integration \cite{butcher}. We can see that the solution obtained using the overset grid is as accurate as the solution obtained without any overset.
\\
\\
\begin{table}
\centering
\resizebox{\textwidth}{!}{%
\begin{tabular}{|c|c|c|c|c|c|}
\hline
\multirow{2}{*}{} &  & \multicolumn{2}{|c|}{DG w/o overset} & \multicolumn{2}{|c|}{DG with overset} \\ \cline{2-6} 
 & Grid size & $L_{2}$ error & Order & $L_{2}$ error & Order  \\ \hline
 \multirow{4}{*}{$\mathbf{P}^{1}$} & 1/20 & 3.215E-03 &  & 3.217E-03 &  \\ \cline{2-6}
  & 1/40 & 7.294E-04 & 2.14 & 7.299E-04 & 2.14 \\ \cline{2-6}
  & 1/80 & 1.725E-04 & 2.08 & 1.726E-04 & 2.08 \\ \cline{2-6}
  & 1/160 & 4.137E-05 & 2.06 & 4.139E-05 & 2.06 \\ \hline
  \multirow{4}{*}{$\mathbf{P}^{2}$} & 1/20 & 2.231E-05 &  & 2.232E-05 &  \\ \cline{2-6}
  & 1/40 & 2.512E-06 & 3.15 & 2.512E-06 & 3.15 \\ \cline{2-6}
  & 1/80 & 2.829E-07 & 3.15 & 2.829E-07 & 3.15 \\ \cline{2-6}
  & 1/160 & 3.187E-08 & 3.15 & 3.187E-08 & 3.15 \\ \hline
  \multirow{4}{*}{$\mathbf{P}^{3}$} & 1/20 & 2.768E-07 &  & 2.771E-07 &  \\ \cline{2-6}
  & 1/40 & 1.376E-08 & 4.33 & 1.387E-08 & 4.32 \\ \cline{2-6}
  & 1/80 & 6.747E-10 & 4.35 & 6.849E-10 & 4.34 \\ \cline{2-6}
  & 1/160 & 3.378E-11 & 4.32 & 3.453E-11 & 4.31 \\ \hline
  \multirow{4}{*}{$\mathbf{P}^{4}$} & 1/20 & 5.890E-09 &  & 5.890E-09 &  \\ \cline{2-6}
  & 1/40 & 1.706E-10 & 5.11 & 1.717E-10 & 5.10 \\ \cline{2-6}
  & 1/80 & 5.294E-12 & 5.01 & 5.366E-12 & 5.00 \\ \cline{2-6}
  & 1/160 & 1.678E-13 & 4.98 & 1.712E-13 & 4.97 \\ \hline
\end{tabular}}
\caption{Validation of overset grid solver using 2D Euler equations for the Isentropic Vortex problem with periodic boundary conditions, $t=10.0$, Uniform mesh with and without overset, $L_{2}$ error for density with $\mathbf{P}^{1}$, $\mathbf{P}^{2}$, $\mathbf{P}^{3}$ and $\mathbf{P}^{4}$ based DG}
\label{table:1}
\end{table}

\section{Results}\label{sec:results}

\begin{figure}[htbp]
  \begin{center}
    \scalebox{0.6}{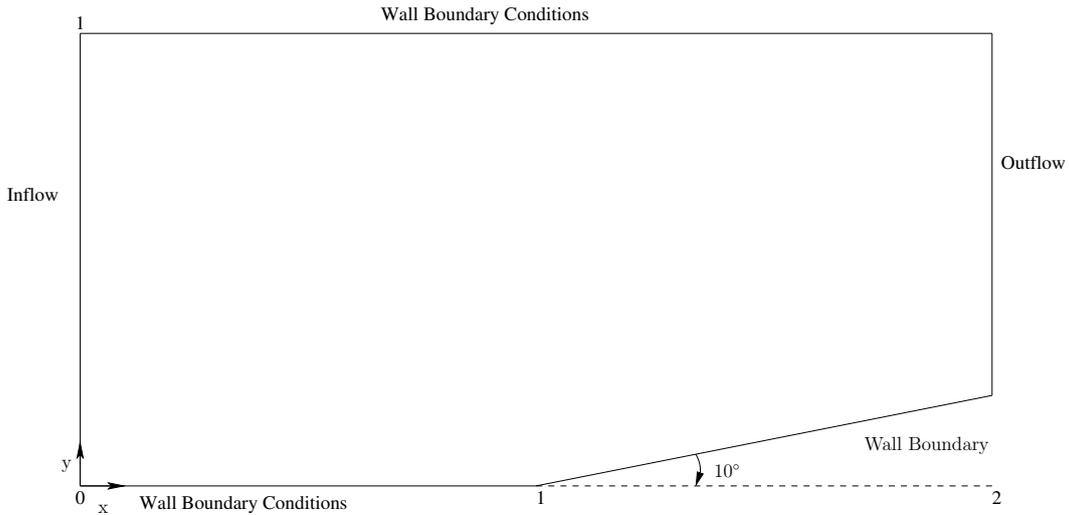}
  \end{center}
  \caption{Computational domain for supersonic flow over a $10^{\circ}$ ramp}
  \label{fig:RampDomain}
\end{figure}

\begin{figure}[htbp]
  \centering
  \subfloat[Density variation for flow over a $10^{\circ}$ ramp for inflow Mach number 3.0]{\label{fig1:rampM3p0}\includegraphics[width=0.9\textwidth]{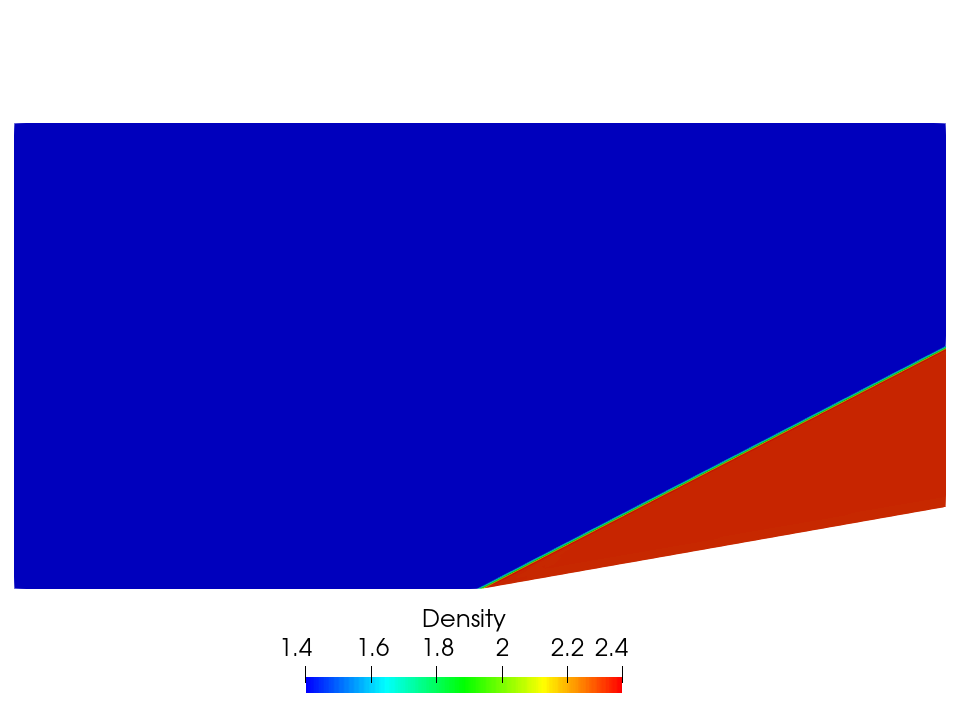}}\hfill
  \subfloat[Zoomed-In view near the base of the ramp showing the grid lines which are colored according to the density values at the cell vertices]{\label{fig2:rampM3p0Zoom}\includegraphics[width=0.9\textwidth]{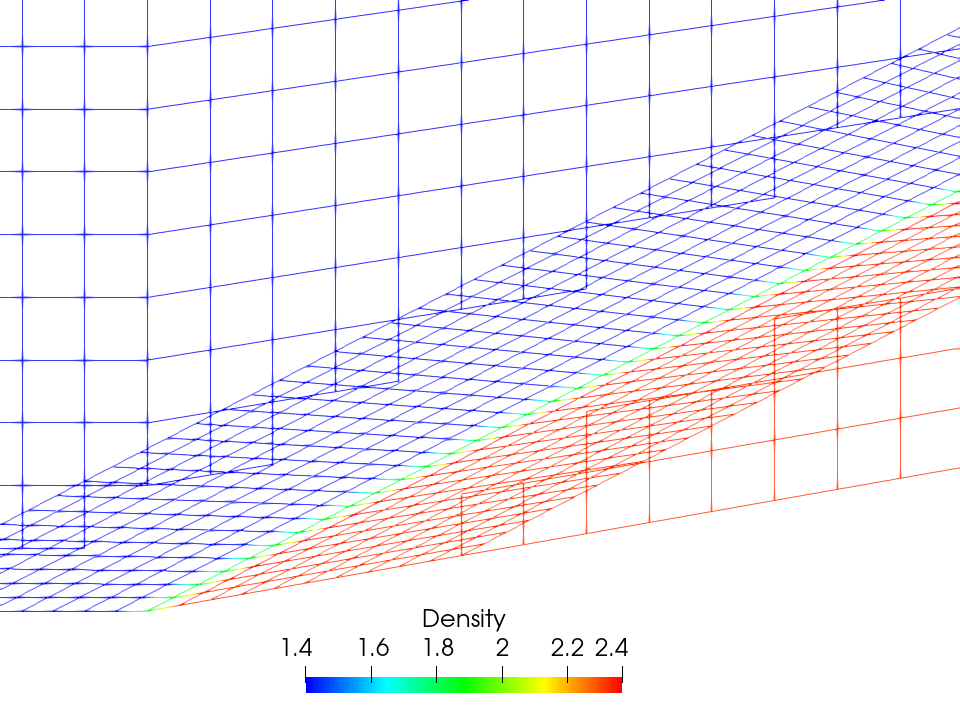}}\hfill
  \caption{Density solution for flow over a $10^{\circ}$ ramp for inflow Mach number 3.0 with $\mathbf{P}^{4}$ based DGM}
  \label{fig:rampM3p0}
\end{figure}

\begin{figure}[htbp]
  \centering
  \subfloat[Density variation for flow over a $10^{\circ}$ ramp for inflow Mach number 4.0]{\label{fig1:rampM4p0}\includegraphics[width=0.9\textwidth]{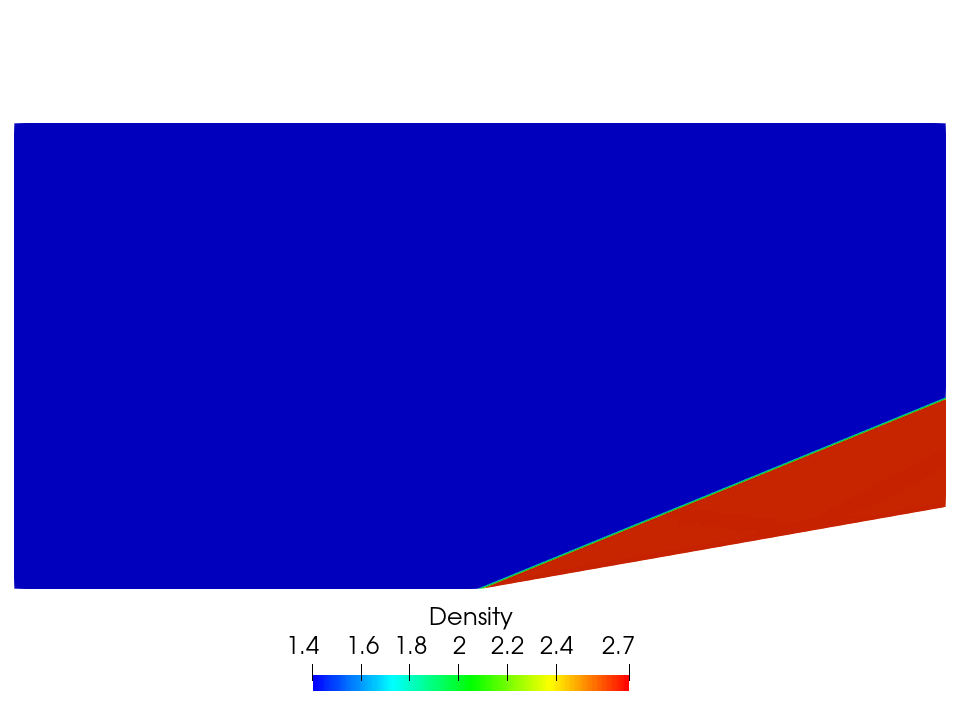}}\hfill
  \subfloat[Zoomed-In view near the base of the ramp showing the grid lines which are colored according to the density values at the cell vertices]{\label{fig2:rampM4p0Zoom}\includegraphics[width=0.9\textwidth]{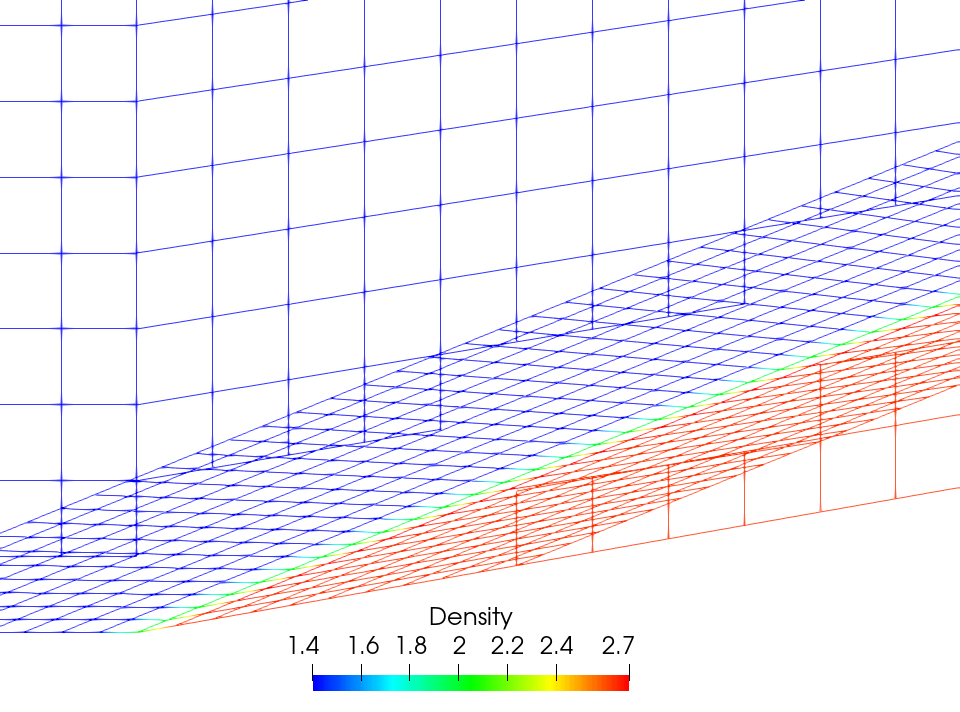}}\hfill
  \caption{Density solution for flow over a $10^{\circ}$ ramp for inflow Mach number 4.0 with $\mathbf{P}^{4}$ based DGM}
  \label{fig:rampM4p0}
\end{figure}

\noindent \textbf{Example 1: Supersonic Flow over a $10^{\circ}$ ramp:} We solve the two dimensional Euler equations as given by \eqref{2dEulerEquations} for supersonic flow over a $10^{\circ}$ ramp. The solution domain is as shown in Figure \ref{fig:RampDomain}. We use an overset grid at the location of the shock which is located using the troubled cells as explained in Section \ref{sec:oversetExpl}. In this particular case, we know that the shock is straight and we keep an appropriate overset grid which covers all the troubled cells. The converged density solutions using $\mathbf{P}^{4}$ based DGM for Mach numbers $3.0$ and $4.0$ are shown in Figures \ref{fig1:rampM3p0} and \ref{fig1:rampM4p0} respectively. A zoomed in view of the solutions near the base of the ramp are shown in Figures \ref{fig2:rampM3p0Zoom} and \ref{fig2:rampM4p0Zoom} showing the grid lines which are colored according to the density values at the cell vertices to show the overset grids and the captured shock. We can clearly see that the shock is aligned to a grid line.
\\
\\

\begin{figure}[htbp]
  \centering
  \subfloat[Density variation for shock reflecting off a flat plate for inflow Mach number 3.0]{\label{fig1:RRM3p0}\includegraphics[width=0.9\textwidth]{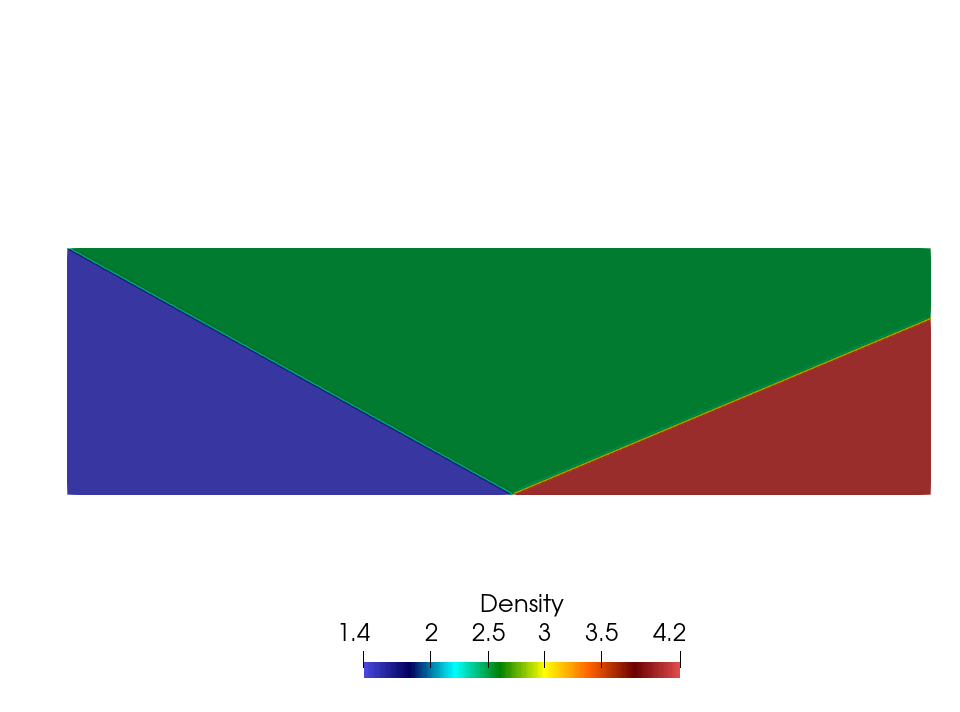}}\hfill
  \subfloat[Zoomed-In view near the impingement point of the shock on the flat plate showing the grid lines which are colored according to the density values at the cell vertices]{\label{fig2:RRM3p0Zoom}\includegraphics[width=0.9\textwidth]{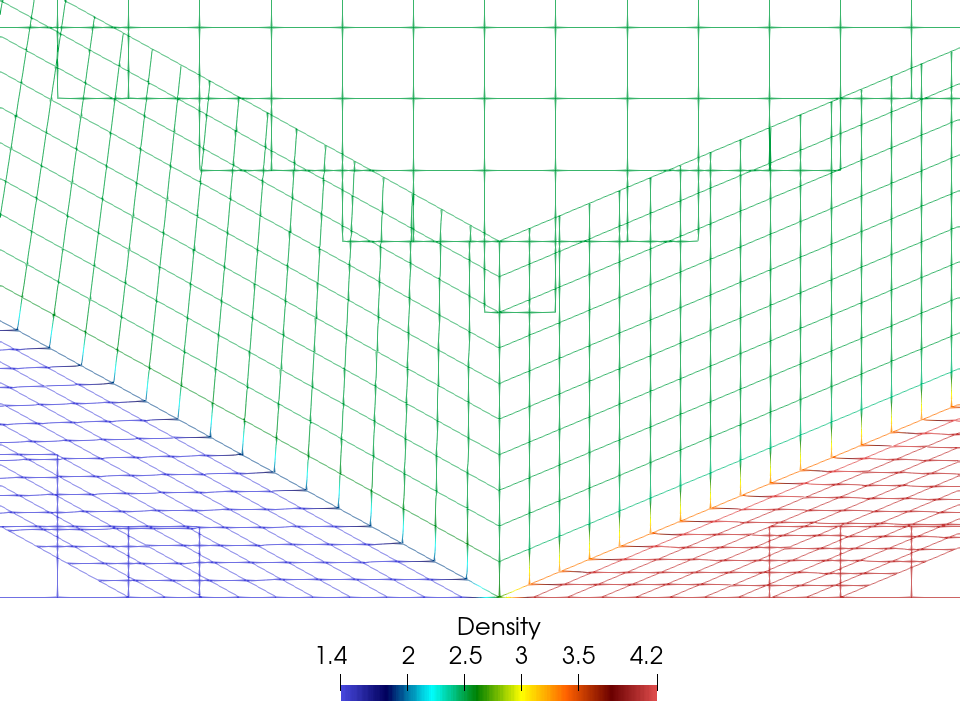}}\hfill
  \caption{Density solution for shock reflecting off a flat plate for inflow Mach number 3.0 with $\mathbf{P}^{4}$ based DGM}
  \label{fig:RRM3p0}
\end{figure}

\begin{figure}[htbp]
  \centering
  \subfloat[Density variation for shock reflecting off a flat plate for inflow Mach number 4.0]{\label{fig1:RRM4p0}\includegraphics[width=0.9\textwidth]{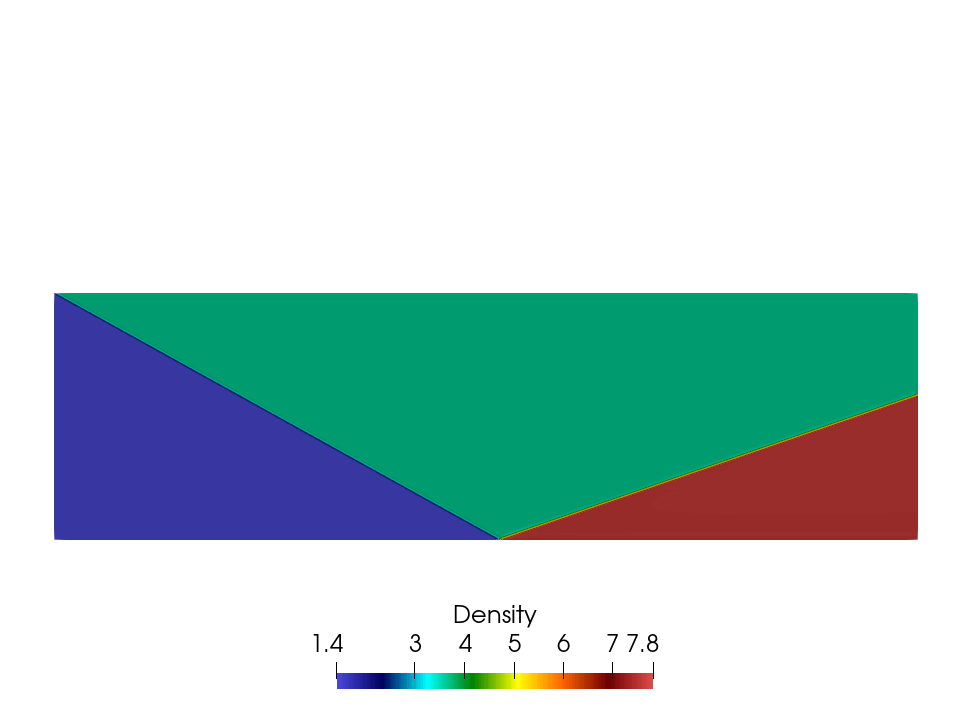}}\hfill
  \subfloat[Zoomed-In view near the impingement point of the shock on the flat plate showing the grid lines which are colored according to the density values at the cell vertices]{\label{fig2:RRM4p0Zoom}\includegraphics[width=0.9\textwidth]{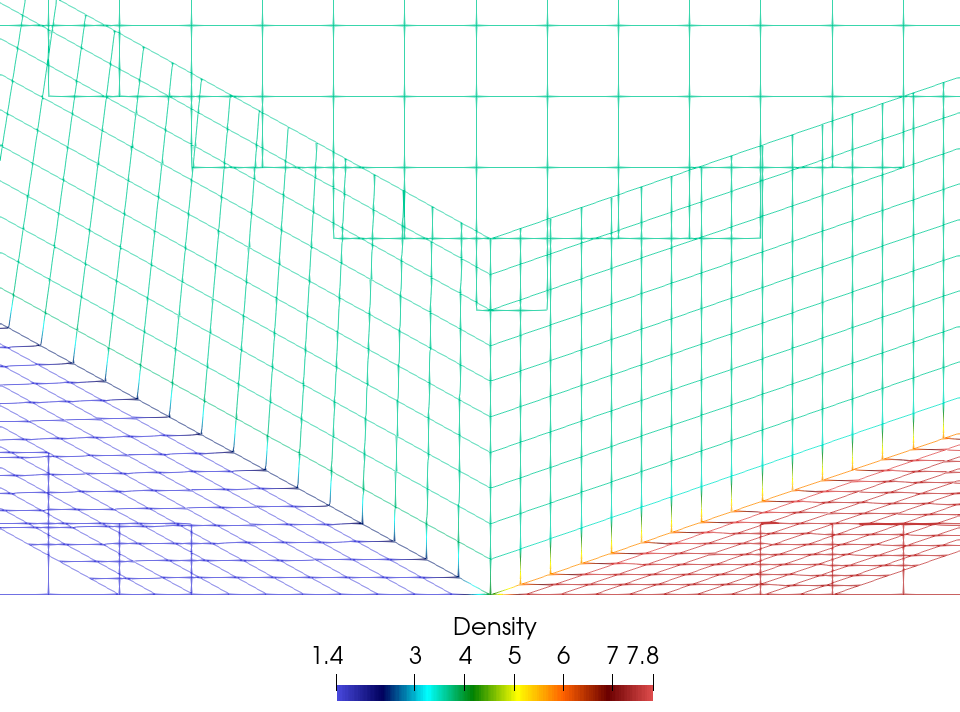}}\hfill
  \caption{Density solution for shock reflecting off a flat plate for inflow Mach number 4.0 with $\mathbf{P}^{4}$ based DGM}
  \label{fig:RRM4p0}
\end{figure}

\noindent \textbf{Example 2: Shock reflecting off a flat plate:} As another test problem for the two-dimensional case with a straight shock, we consider a flow where a shock reflects off a flat plate with inflow Mach numbers $3.0$ and $4.0$. We solve the two-dimensional Euler equations \eqref{2dEulerEquations} in the computational domain $[0,3.5]\times[0,1]$. For the left boundary, we impose the inflow conditions. For the right boundary, we use supersonic exit boundary conditions. For the bottom boundary, we impose wall boundary conditions and for the top boundary, we use the post shock conditions for inflow Mach number and an oblique shock of shock angle $29^{\circ}$. Again, from the troubled cells and knowing that the shock is straight, we construct an overset grid which covers all the troubled cells. The converged density solutions using $\mathbf{P}^{4}$ based DGM for inflow Mach numbers $3.0$ and $4.0$ are shown in Figures \ref{fig1:RRM3p0} and \ref{fig1:RRM4p0} respectively.  A zoomed in view of the solutions near the impingement point of the shock on the flat plate are shown in Figures \ref{fig2:RRM3p0Zoom} and \ref{fig2:RRM4p0Zoom} showing the grid lines which are colored according to the density values at the cell vertices to show the overset grids and the captured shock. Again, we can clearly see that the shock is aligned to a grid line.
\\
\\

\begin{figure}[htbp]
  \begin{center}
    \scalebox{0.6}{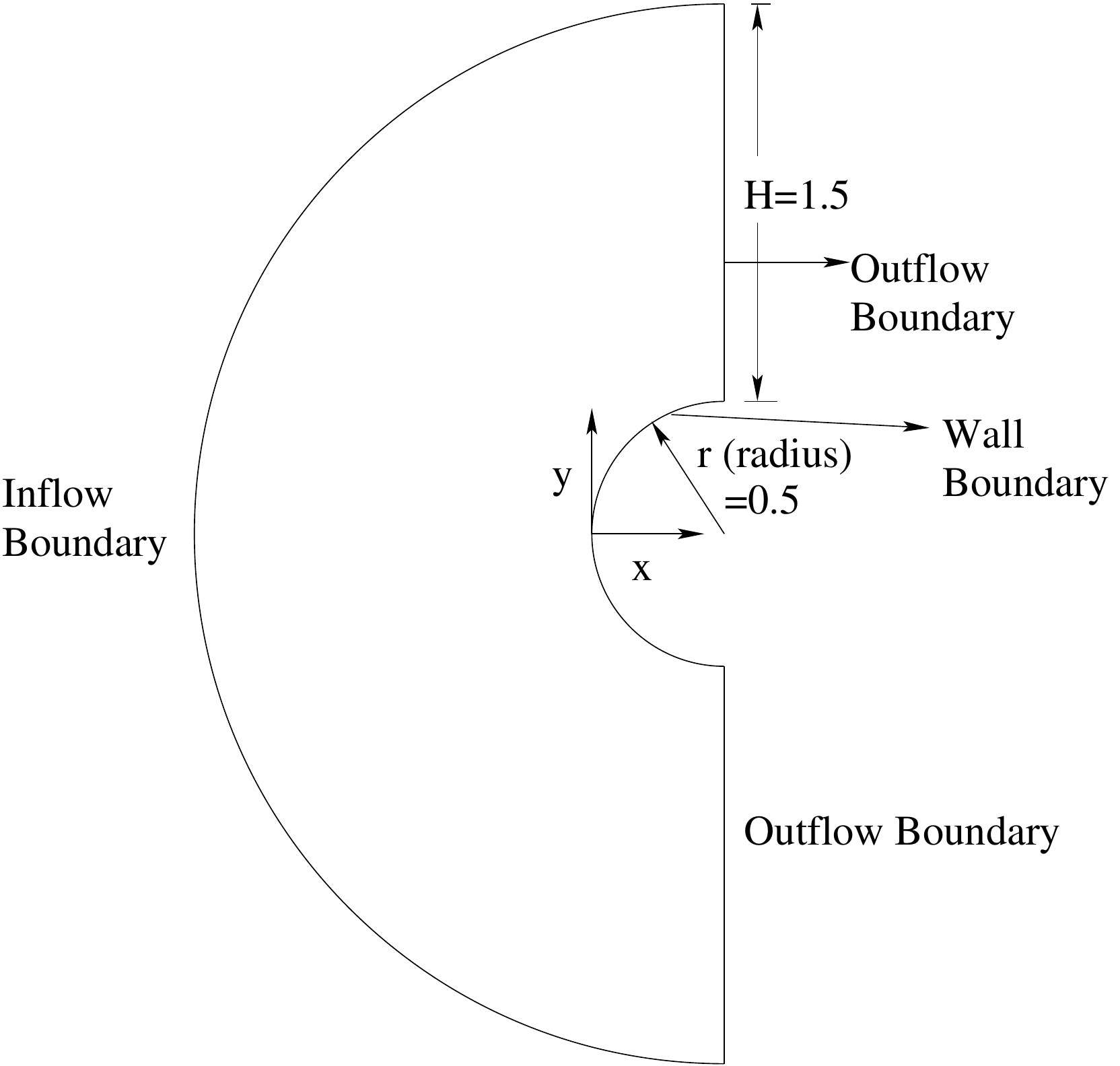}
  \end{center}
  \caption{Computational domain for Supersonic flow over a circular cylinder}
  \label{fig:CylDomain}
\end{figure}

\begin{figure}[htbp]
  \centering
  \subfloat[Mach number variation for flow over a circular cylinder for inflow Mach number 3.0]{\label{fig1:CircCylM3p0}\includegraphics[width=0.9\textwidth]{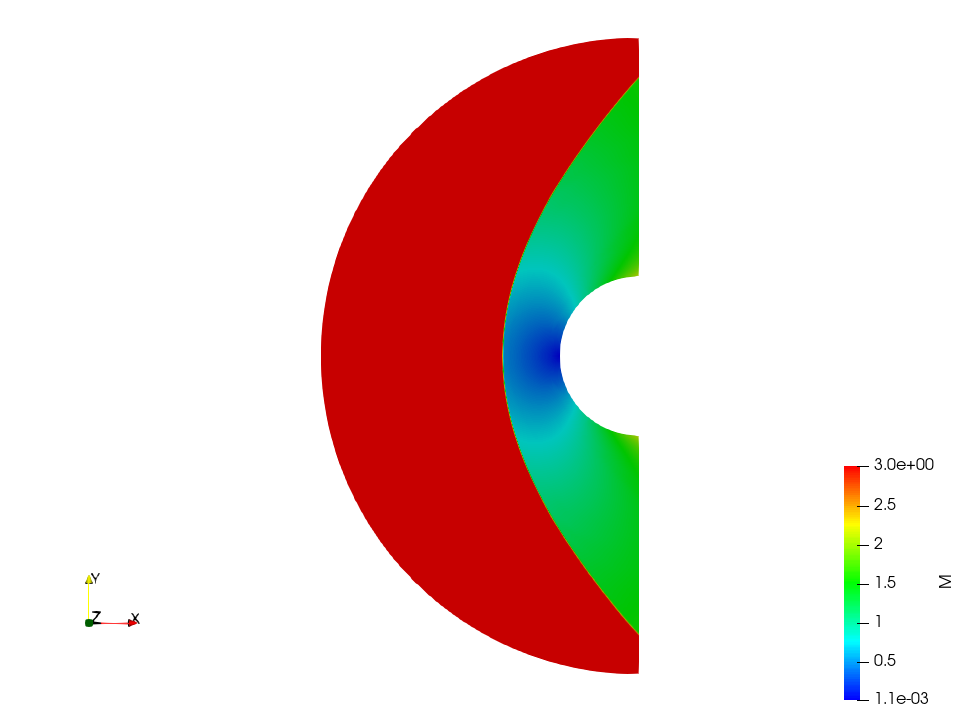}}\hfill
  \subfloat[Zoomed-In view near the location of the shock showing the grid lines which are colored according to the Mach number values at the cell vertices according to the color scheme shown in (a)]{\label{fig2:CircCylM3p0Zoom}\includegraphics[width=0.9\textwidth]{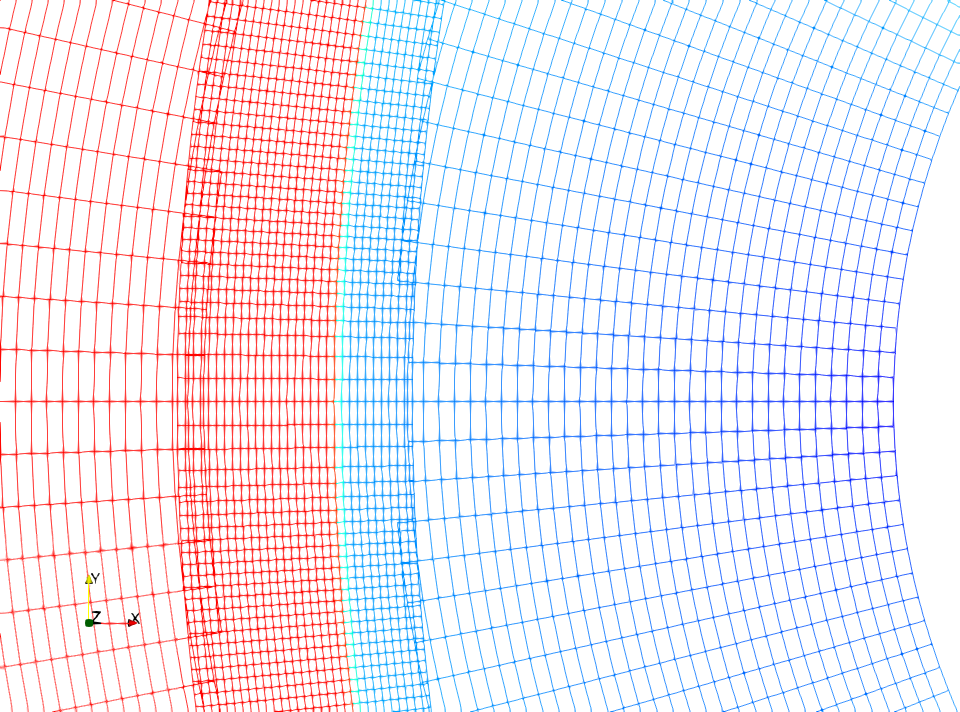}}\hfill
  \caption{Mach number solution for flow over a circular cylinder of radius 0.5 for inflow Mach number 3.0 with $\mathbf{P}^{4}$ based DGM}
  \label{fig:CircCylM3p0}
\end{figure}

\begin{figure}[htbp]
  \centering
  \subfloat[Mach number variation for flow over a circular cylinder for inflow Mach number 4.0]{\label{fig1:CircCylM4p0}\includegraphics[width=0.9\textwidth]{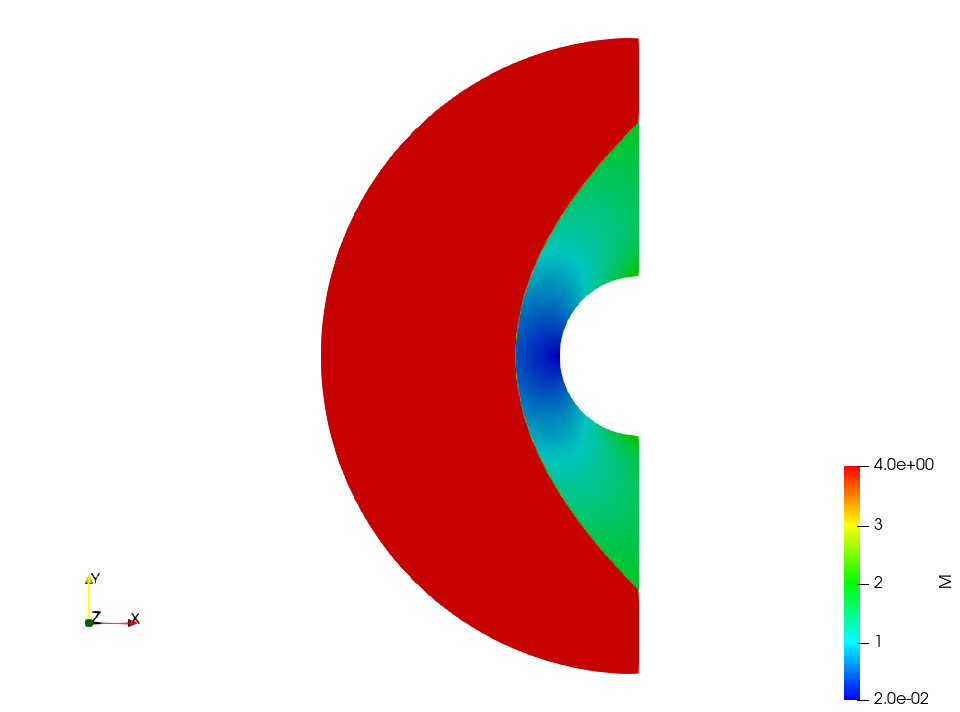}}\hfill
  \subfloat[Zoomed-In view near the location of the shock showing the grid lines which are colored according to the Mach number values at the cell vertices according to the color scheme shown in (a)]{\label{fig2:CircCylM4p0Zoom}\includegraphics[width=0.9\textwidth]{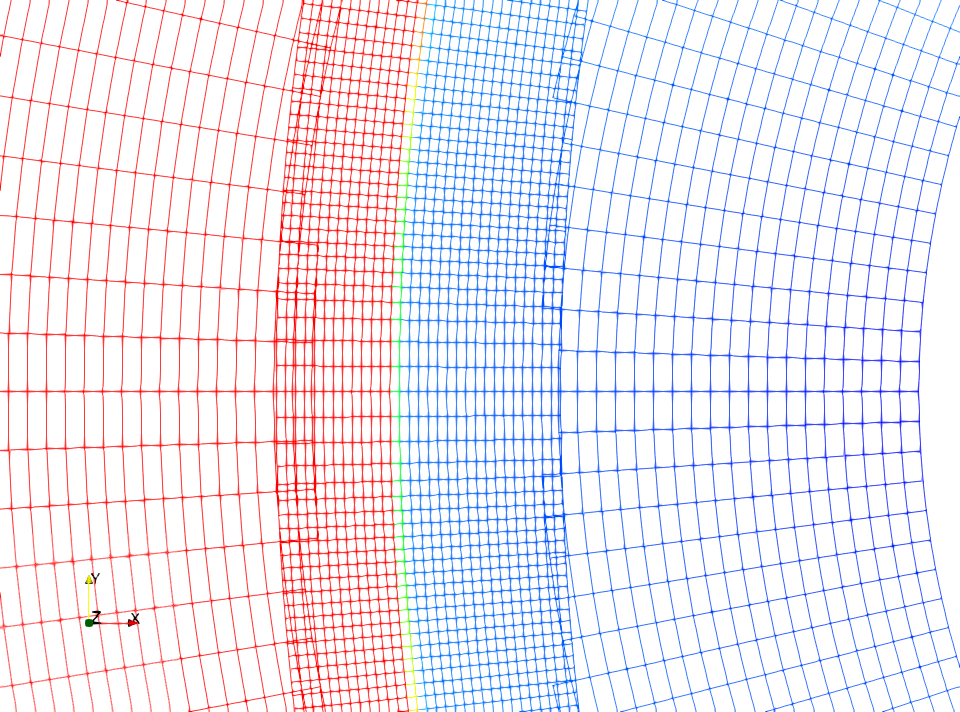}}\hfill
  \caption{Mach number solution for flow over a circular cylinder of radius 0.5 for inflow Mach number 4.0 with $\mathbf{P}^{4}$ based DGM}
  \label{fig:CircCylM4p0}
\end{figure}

\begin{figure}[htbp]
\begin{center}
\includegraphics[scale=0.75]{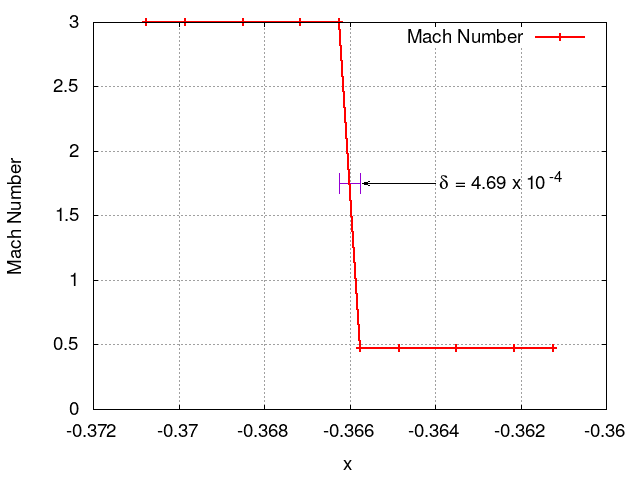}
\caption{Variation of Mach number as a function of $x$ along the center line ($y$ = 0) for the two cells near the shock showing the solution at the Gauss Legendre quadrature points for flow over a circular cylinder with inflow Mach number 3.0 with a non-dimensional shock thickness $\delta = 4.69\times 10^{-4}$ using $\mathbf{P}^{4}$ based DGM.}
\label{fig:CircCylM3p0CenterMVariation}
\end{center}
\end{figure}

\begin{figure}[htbp]
\begin{center}
\includegraphics[scale=0.75]{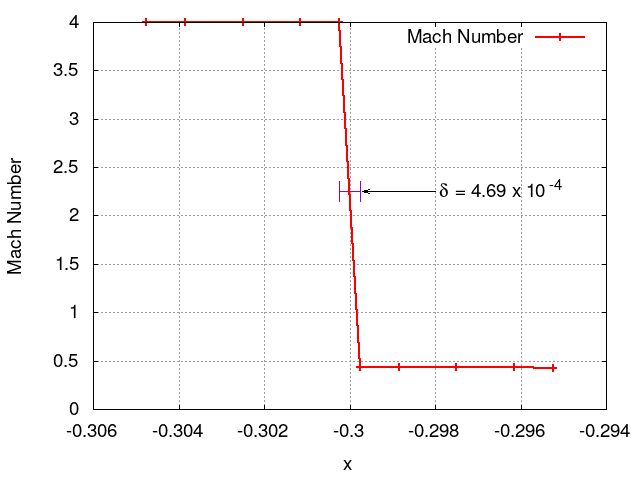}
\caption{Variation of Mach number as a function of $x$ along the center line ($y$ = 0) for the two cells near the shock showing the solution at the Gauss Legendre quadrature points for flow over a circular cylinder with inflow Mach number 4.0 with a non-dimensional shock thickness $\delta = 4.69\times 10^{-4}$ using $\mathbf{P}^{4}$ based DGM.}
\label{fig:CircCylM4p0CenterMVariation}
\end{center}
\end{figure}

\noindent \textbf{Example 3: Supersonic flow over a circular cylinder:} As a test problem with a curved shock, we consider the supersonic flow over a circular cylinder of radius $0.5$. We solve the two-dimensional Euler equations \eqref{2dEulerEquations} in the computational domain which is shown in Figure \ref{fig:CylDomain}. The troubled cells obtained using a coarse grid of $150$ (circumferential) and $\times 100$ (wall-normal) elements for inflow Mach number $3.0$ are shown in Figure \ref{fig2:CircCylM3p0troubledCells}. We use the troubled cell data and take the first troubled cell encountered by the flow and fit a curve using a second order b-spline following the procedure outlined in Section \ref{sec:oversetExpl}. We construct an overset grid using this b-spline and curves parallel to this b-spline as shown in Figure \ref{fig2:CircCylM3p0OversetGrid1}. We reduce the width of the overset grid even further using the troubled cells as explained in Section \ref{sec:oversetExpl}. With this overset grid, we obtained the solution for inflow Mach numbers $3.0$, $3.2$, $3.4$, $3.6$, $3.6$ and $4.0$. The Mach number solution obtained using $\mathbf{P}^{4}$ based DGM where the shock is aligned with a grid line is shown for Mach numbers $3.0$  and $4.0$ in Figures \ref{fig1:CircCylM3p0}
and \ref{fig1:CircCylM4p0} respectively. A zoomed in view of the solutions near the shock are shown in Figures \ref{fig2:CircCylM3p0Zoom}, 
and \ref{fig2:CircCylM4p0Zoom} showing the grid lines colored according to the Mach number values at the cell vertices. Here, we can see that the shock aligns with a grid line. We have also shown the variation of Mach number along the center line ($y$=0) near the location of the shock for inflow Mach numbers $3.0$ and $4.0$ in Figures \ref{fig:CircCylM3p0CenterMVariation} and \ref{fig:CircCylM4p0CenterMVariation} respectively. From these figures, we can clearly see that the shock has been captured quite sharply with a non-dimensional shock thickness of $\delta = 4.69\times 10^{-4}$. We note that this thickness depends on the order of the scheme. If we increase the order, this will reduce further. To verify the solution obtained, we used the method proposed in \cite{sc} to calculate the shock offset distance $x_{O}$ and the location of the sonic point on the body $(x_{SB},y_{SB})$ using the coordinate system shown in Figure \ref{fig:CylDomain}. These values are compared with the values obtained from our simulations in Table \ref{table:2}. We can see that the shock offset distance and the location of sonic point on the body obtained from our simulations agree quite well with the analytically calculated value.
\\
\\

\begin{table}
\centering
\resizebox{\textwidth}{!}{%
\begin{tabular}{|c|c|c|c|c|c|c|}
\hline
\multirow{2}{*}{} & \multicolumn{3}{|c|}{Analytical \cite{sc}} & \multicolumn{3}{|c|}{Computational} \\ \cline{2-7} 
 M & -$x_{O}$ & $x_{SB}$ & $y_{SB}$ & -$x_{O}$ & $x_{SB}$ & $y_{SB}$  \\ \hline
 3.0 & 0.367 & 0.152 & 0.359  & 0.366 & 0.151 & 0.360 \\ \hline
 3.2 & 0.347 & 0.150 & 0.357  & 0.347 & 0.151 & 0.358 \\ \hline
3.4 & 0.331 & 0.148 & 0.355  & 0.330 & 0.148 & 0.356 \\ \hline
  3.6 & 0.318 & 0.147 & 0.354  & 0.320 & 0.147 & 0.354 \\ \hline
  3.8 & 0.308 & 0.146 & 0.353  & 0.308 & 0.146 & 0.352 \\ \hline
  4.0 & 0.299 & 0.145 & 0.352  & 0.300 & 0.144 & 0.352 \\ \hline
\end{tabular}}
\caption{Shock Offset Distance (-$x_{O}$) and the location of sonic point $(x_{SB},y_{SB})$ on the body for flow over a circular cylinder for various Mach numbers}
\label{table:2}
\end{table}



\section{Conclusion:}\label{sec:conc}

\noindent We have developed a new procedure for shock capturing using overset grids with discontinuous Galerkin method. We run the solver on a coarse grid and use the troubled cell data to construct an overset grid which is aligned to the shock and also encompasses all the troubled cells. This is done using b-spline interpolation for a curved shock. We run the solver again using the coarse grid solution as the initial condition and using the limiter only in the overset grid. We also use a high resolution numerical flux in the overset grid. This procedure gives us a shock aligned with a grid line capturing it sharply with a very small non-dimensional shock thickness. We have used this method to capture the straight shocks for supersonic flow over a $10^{\circ}$ ramp and shock reflecting over a flat plate. We have also captured the shock for flow over a circular cylinder. This method will work for any flow containing a discontinuity to capture the discontinuity aligned to a grid line.


\bibliographystyle{ieeetr}
\bibliography{references}

\begin{thebibliography}{10}

\bibitem{salas1}
M.~D. Salas, {\em {A Shock-Fitting Primer}}.
\newblock {CRC Applied Methematics and Nonlinear Science, Chapman and Hall, 1st
  edition}, 2009.

\bibitem{bpm}
A.~Bonfiglioli, R.~Paciorri, and A.~D. Mascio, ``{The Role of Mesh Generation,
  Adaptation, and Refinement on the Computation of Flows Featuring Strong
  Shocks}.,'' {\em Modelling and Simulation in Engineering}, 631276, 2012.

\bibitem{qs2}
J.~Qiu and C.-W. Shu, ``A comparison of troubled-cell indicators for
  {Runge-Kutta} discontinuous {Galerkin} methods using weighted essentially
  nonoscillatory limiters.,'' {\em SIAM J. Sci. Comput.}, vol.~27,
  pp.~995--1013, 2005.

\bibitem{qs1}
J.~Qiu and C.-W. Shu, ``{Runge-Kutta} discontinuous {Galerkin} method using
  {WENO} limiters.,'' {\em SIAM J. Sci. Comput.}, vol.~26, pp.~907--929, 2005.

\bibitem{pp}
P.-O. Persson and J.~Peraire, ``{Sub-cell} shock capturing for discontinuous
  {Galerkin} methods.,'' {\em AIAA 2006-112}, 2006.

\bibitem{yh}
J.~Yu and J.~Hesthaven, ``A comparative study of shock capturing models for the
  discontinuous galerkin method.,'' 2017.

\bibitem{gq}
P.~Giri and J.~Qiu, ``{A high-order Runge-Kutta discontinuous Galerkin method
  with a subcell limiter on adaptive unstructured grids for two-dimensional
  compressible inviscid flows},'' {\em International Journal for Numerical
  Methods in Fluids}, vol.~91, no.~8, pp.~367--394, 2019.

\bibitem{bsd}
J.~Benek, J.~Steger, and F.~Dougherty, ``A {Flexible Grid Embedding Technique
  with Application to the Euler Equations}.,'' {\em AIAA Paper 1983-1944},
  1983.

\bibitem{butcher}
J.~C. Butcher, {\em {Numerical Methods for Ordinary Differential Equations}}.
\newblock {John Wiley and Sons}, 2016.

\bibitem{lpr}
W.~Li, J.~Pan, and Y.-X. Ren, ``The discontinuous {Galerkin spectral element
  methods for compressible flows on two-dimensional mixed grids}.,'' {\em
  Journal of Computational Physics}, vol.~364, pp.~314--346, 2018.

\bibitem{gbot}
M.~Galbraith, J.~Benek, P.~Orkwis, and M.~Turner, ``A {Discontinuous Galerkin
  Chimera scheme}.,'' {\em Computers and Fluids}, vol.~98, pp.~27--53, 2014.

\bibitem{ks3}
K.~Kitamura and E.~Shima, ``Towards {shock-stable} and accurate hypersonic
  heating computations: {A new pressure flux for AUSM-family schemes}.,'' {\em
  Journal of Computational Physics}, vol.~245, pp.~62--83, 2013.

\bibitem{shu1}
C.-W. Shu, ``Essentially non-oscillatory and weighted essentially
  non-oscillatory schemes for hyperbolic conservation laws.,'' {\em Lecture
  Notes in Mathematics, Springer}, vol.~1697, pp.~325--432, 1998.

\bibitem{sc}
J.~Sinclair and X.~Cui, ``A theoretical approximation of the shock standoff
  distance for supersonic flows around a circular cylinder.,'' {\em Phys.
  Fluids}, vol.~29, 026102, 2017.

\end{thebibliography}

\end{document}